%% file: splice.tex
\documentstyle[11pt, sec, psfig, diagram]{article}

\newtheorem{theorem}{Theorem}[section]
\newtheorem{lemma}[theorem]{Lemma}
\newtheorem{proposition}[theorem]{Proposition}
\newtheorem{corollary}[theorem]{Corollary}
\newtheorem{remark}[theorem]{Remark}
\newtheorem{definition}[theorem]{Definition}

\def\ra{\rightarrow}

\def\be{\begin{equation}}
\def\ee{\end{equation}}

\input{dep.tex}

\begin{document}

\title{On the Cappell-Lee-Miller glueing theorem}

\author{Liviu I.Nicolaescu\thanks{{\bf Current address}: Dept.of Math., McMaster University, Hamilton, Ontario,  L8S 4K1, Canada; nicolaes@icarus.math.mcmaster.ca}}

\date{Version 3: March, 1998}

\maketitle

\begin{abstract} 
We formulate a more  conceptual interpretation  of the  Cappell-Lee-Miller  glueing/splitting theorem  in terms  of asymptotic maps and asymptotic exact sequences.   Additionally, we show this gluing result is equivalent to a Mayer-Vietoris-type long exact sequence.  We also present   applications   to   eigenvalue estimates, approximation of obstruction bundles   and glueing of determinant line bundles arising  frequently in gauge theory.

 All these results are true  in a slightly more general context than in \cite{CLM}. We work   with operators which    differ  from  translation invariant ones by  exponentially decaying terms. 
\end{abstract}

\addcontentsline{toc}{section}{Introduction}

\begin{center}
{\bf Introduction}
\end{center}

\bigskip

\input{spli0.tex}

\tableofcontents

\section{First order elliptic equations on manifolds with cylindrical ends}
\setcounter{equation}{0}
\input{spli1.tex}

\section{Proof of the Main Theorem}
\setcounter{equation}{0}
\input{spli2.tex}

\section{A Mayer-Vietoris interpretation}
\input{spli3.tex}

\section{Applications}
\setcounter{equation}{0}
\input{spli4.tex}

\appendix
\section{Some technical proofs}
\setcounter{equation}{0}
\input{splia.tex}

\end{document}

%% file: dep.tex
\newfont{\Bbb}{msbm10 scaled \magstep1}

\newcommand\bR{\hbox{\Bbb R}}

\newcommand\bZ{\hbox{\Bbb Z}}

\newfont{\es}{eusm10 scaled \magstep1}
\newfont{\ses}{eufm8 scaled \magstep1}
\newfont{\gt}{eufb10 scaled \magstep1}
\newfont{\sg}{eufb8 scaled \magstep1}
\newfont{\goth}{eufb10 scaled \magstep2}

\newcommand{\dir}{\hbox{\es D}}
\newcommand{\cdir}{ {\not\!\!{\hbox{\es D}}}}

\def\e{\cal E}

\def\k{\cal K}

\def\ra{\rightarrow}

\newcommand{\Lra}{{\longrightarrow}}

\def\be{\begin{equation}}
\def\ee{\end{equation}}

\def\lan{\langle}
\def\ran{\rangle}

\newcommand{\ve}{{\varepsilon}}



%% file: spli0.tex
  Consider the following set-up. We are given  two oriented, Riemann manifolds  $M_i(\infty)$, $i=1,2$. $M_1(\infty)$ has a   (metrically) cylindrical end ${\bR}_+\times N$  while $M_2(\infty)$ has a cylindrical end ${\bR}_-\times N$. Here  $N$ is a closed, compact oriented  Riemann manifold, not necessarily connected.   $\hat{E}_i\ra M_i(\infty)$ are  bundles equipped with inner products along their fibers and  $\hat{D_i}:C^\infty(\hat{E})\ra C^\infty(\hat{E}_i)$ are     self-adjoint Dirac-type operators   which along the necks have the form 
\[
\hat{D_i}= J(\nabla_t- D)
\]
where  $D$ is a selfadjoint Dirac-type operator  operator on the $E:= E_1\!\mid_{N}\cong \hat{E}_2\!\mid_N$ and $J$ denotes the Clifford multiplication by $dt$. We want to emphasize  that $D$ is independent of the longitudinal coordinate $t$ along the necks.   Consider now two  smooth self-adjoint endomorphisms  $\hat{B}_i$ of $\hat{E}_i$  and along the necks  set 
\[
B_i(t):= \hat{B}_i\!\mid_{t\times N},\;\; A_i(t):= JB_i(t).
\]
We assume the  following  about  $A_i(t)$.

\medskip

\begin{figure}
\centerline{\psfig{figure=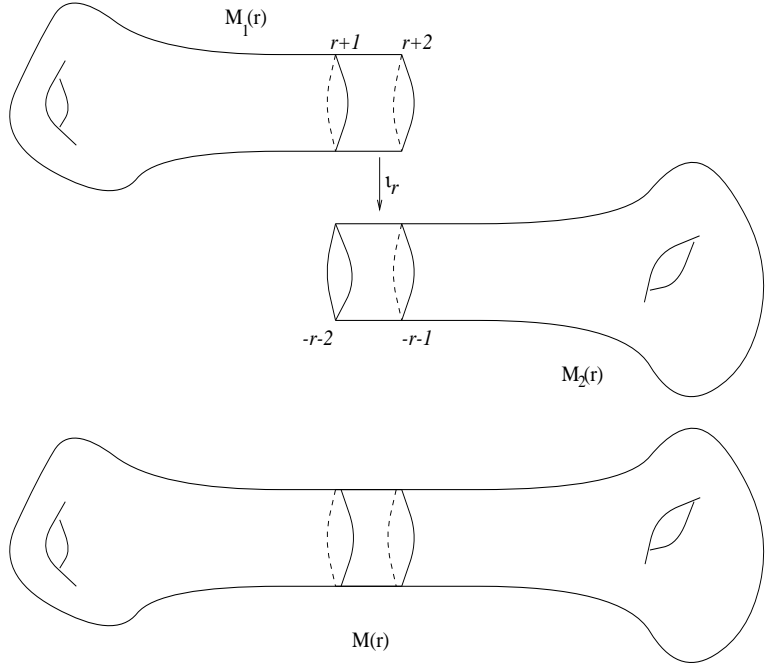,height= 3.5in,width=4in}}
\caption{\sl Glueing two manifolds with cylindrical ends.}
\label{fig: spli2}
\end{figure}

\noindent $\bullet$ $A_i(t)$ anti-commutes with  $J$
\[
\{J, A_i(t)\} = JA_i(t)+A_i(t) J=0.
\]
\noindent $\bullet$ There exist $C, \lambda >0$  such that 
\be
\sup\{ |\hat{A}_i(x)| \; ;\;x\in [t,t+1]\times N\} \leq C \exp (-\lambda |t|).
\label{eq: decay1}
\ee
Consider a smooth, decreasing, cut-off function $\eta: {\bR}_+ \ra [0,1]$ such that  
\[
\eta(t) \equiv 1, \;\; t\in [0,1/4]
\]
\[
\eta(t) \equiv 0,\;\; t\geq 3/4
\]
and 
\[
\left|\frac{d\eta}{dt}\right| \leq 4,\;\;\forall t\geq 0.
\]
For each $r>0$ set $\eta_r(t):=\eta(t-r)$. Now extend $\eta_r$ by symmetry to a function on ${\bR}$ still denoted by $\eta_r$. We can regard $\eta_r(t)$, $t\geq 0$,  as a smooth function on $M_1(\infty)$  and $\eta_r(t)$, $t\leq 0$  as a smooth function on $M_2(\infty)$  so we can form
\[
\hat{D}_{i,r}:= \hat{D}_i + \eta_r\hat{B}_i.
\]
Also define
\[
\hat{D}_{i,\infty} :=\hat{D}_i+\hat{B}_i.
\]
Note that  via the  obvious linear increasing  diffeomorphism 
\[
\imath_r:[r+1,r+2]\ra [-r-2,-r-1]
\]
 we have identifications
\[
D_{1,r}\!\mid_{C(r+1)}  = D_{2,r}\!\mid_{C(-r-2)}, \;\;C(t):=[t,t+1]\times N.
\]
Denote by $M_1(r)$ the manifold obtained from $M_1(\infty)$ by chopping off the cylinder $(r+2,\infty)\times N$   and by $M_2(r)$ the manifold  $M_2(\infty)\setminus (-\infty,-r-2)\times N$. We glue $M_1(r)$ to $M_2(r)$ via $\imath_r$ to obtain a closed manifold $M(r)$; see Figure \ref{fig: spli2}.  Similarly the operators $\hat{D}_{i,r}$ can be glued together to produce  a Dirac type operator $\dir_r$ on the  obvious  glued bundle ${\cal E}_r\ra M(r)$.

The operators $\hat{D}_i$ may have additional symmetries. We will be particularly interested in {\em super-symmetric} operators. This means  the bundles $\hat{E}_i$ are equipped  with  orthogonal (unitary) decompositions
\be
\hat{E}_i=\hat{E}_i^+\oplus \hat{E}_i^-
\label{eq: super}
\ee
which determine the chiral operators $\hat{C}_i:= \hat{P}^+_i -\hat{P}_i^-$ where  $\hat{P}_i^\pm$ denotes the orthogonal projection $\hat{E}_i\ra \hat{E}_i^\pm$.   The Dirac operator $\hat{D}_i$ is said to be super-symmetric  if
\be
\{ C_i, \hat{D}_i\}=0.
\label{eq: super1}
\ee
Equivalently, in terms of the splitting (\ref{eq: super})  it has the  block decomposition
\[
\hat{D}_i=\left[
\begin{array}{cc}
0 & \hat{\bf D}_i^*
\\
\hat{\bf D}_i & 0 
\end{array}
\right]
\]
where $\hat{\bf D}_i $ is a first order  elliptic operator $C^\infty(\hat{E}_i^+)\ra C^\infty(\hat{E}_i^-)$.  The condition (\ref{eq: super1}) implies that  for any  1-form $\alpha$ on $\hat{M}_i(\infty)$ the Clifford multiplication by $\alpha$ anti-commutes with $\hat{C}_i$
\be
\{\hat{c}(\hat{\alpha}), \hat{C}_i\}=0.
\label{eq: super2}
\ee
Note that along the neck the operator $\hat{\bf D}_i$ has the form 
\[
\hat{\bf D}_i=G(\nabla_t -{\bf D})
\]
where $G:E^+\ra E^-$ is the bundle isomorphism given by the Clifford multiplication by $dt$ and ${\bf  D}:C^\infty(E_i^+)\ra C^\infty(E_i^+)$  is a self-adjoint, Dirac-type operator.

We will further assume  that the two super-symmetries are compatible along the ``boundary'' $N$ i.e. 
\[
\hat{C}_1\!\mid_N=\hat{C}_2\!\mid_N =: C.  
\]
Thus the bundle $E$ is super-symmetric with chiral operator $C$. The  conditions (\ref{eq: super1}) and  (\ref{eq: super2}) imply that 
\be
[C, D]= CD-DC = 0.
\label{eq: super3}
\ee
In this case we assume the perturbations $\hat{B}_i$ are compatible with the chiral operators in an obvious sense. Clearly  the super-symmetry is transmitted to the glued  bundle ${\cal E}(r)$ and the glued  operator ${\dir}_r$. The kernel ${\cal K}_r$ is naturally a finite dimensional ${\bZ}_2$-graded  space.   In this paper we will address the following  question.

\medskip

\noindent {\bf Main Problem}  {\em Understand the behavior of ${\cal K}_r$ as $r\ra \infty$.}

\medskip

The kernel of an operator is a notoriously  unstable  object  so  it is unrealistic to be able to solve  the Main Problem as stated. We need to ``stabilize'' ${\cal K}_r$ if we  expect to say something of significance.

To formulate the main result we need to introduce  some  additional notions.  We begin with  the notions of asymptotic map and  asymptotic exactness.  An  {\em asymptotic map} is    a sequence $(U_r, V_r,f_r)_{r>0}$ with the following properties

\medskip

(a) There exist Hilbert spaces $H_0$ and $H_1$ such that $U_r$ is a closed subspace of $H_0$ and $V_r$ is a closed subspace of $H_1$, $\forall r>0$.

(b) $f_r$ is a densely defined  linear map $f_r:U_r\ra H_1$ with closed graph and range $R(f_r)$, $\forall r>0$.

(c) $\lim_{r\ra \infty}\hat{\delta}(R(f_r), V_r) =0$ where, following \cite{Kato}, we  set
\[
\hat{\delta}(U, V)= \sup\{ {\rm dist}\, (u, V)\; ;\; u\in U, \;|u|=1\}.
\]
We will  denote  asymptotic maps by $U_r\stackrel{f_r}{\Lra^a}V_r$. There is a super-version of this notion when $U_r$ and $V_r$ are ${\bZ}_2$-graded and  are closed  subspaces in   ${\bZ}_2$-graded Hilbert spaces  such that the natural inclusions are even.

The next result, proved in the Appendix,  explains the motivation behind the above definition.

\begin{lemma}{\rm If 
\[
\hat{\delta}(U, V)<1
\]
then the orthogonal projection $P_V$ onto $V$ induces a  one-to one map $U\ra V$. If additionally
\[
\hat{\delta}(V,U)<1
\]
then $P_V:U\ra V$ is a linear isomorphism.}
\label{lemma: asym}
\end{lemma}

Define the gap between two  closed subspaces $U,V$  in a Hilbert space $H$ by
\[
\delta(U, V)=\max(\hat{\delta}(U,V), \hat{\delta}(V,U)\}.
\]

The sequence of asymptotic maps
\[
U_r\stackrel{f_r}{\Lra^a}V_r\stackrel{g_r}{\Lra^a}W_r,\;\;r\ra \infty
\]
is  said to be  {\em asymptotically exact} if
\[
\lim_{r\ra \infty} \delta( R(f_r), \ker g_r) =0.
\]
We have the following consequence of Lemma \ref{lemma: asym}.

\begin{lemma}{\rm If the sequence
\[
U_r\stackrel{f_r}{\Lra^a}V_r\stackrel{g_r}{\Lra^a}W_r,\;\;r\ra \infty
\]
is asymptotically exact, $P_r$ denotes the orthogonal projection onto  $\ker g_r$ and $Q_r$ the orthogonal projection onto $W_r$ then  there exists $r_0>0$ such that the sequence 
\[
U_r \stackrel{P_r\circ f_r}{\Lra}V_r\stackrel{Q_r\circ g_r}{\Lra}W_r
\]
is exact for  all $r>r_0$.}
\end{lemma}
If  the spaces $H_j$ are ${\bZ}_2$-graded $H_j=H_j^+\oplus H_j^-$ we  say the sequence is {\em super-symmetric} if   the maps $f_r$  and $g_r$ are even i.e. are compatible with the splitting. In this case we get two asymptotically exact sequences
\[
U_r^\pm \ra V^\pm_r\ra W_r^\pm.
\]
Next we need to introduce suitable functional spaces. For brevity we  discuss only   distributions on $M_1(\infty)$. Define the {\em extended} $L^2$-space $L^2_{ex}(\hat{E}_1)$ as the space of sections $\hat{u}\in L^2_{loc}(\hat{E}_1)$ such that  there exists $u_\infty\in L^2(E)$ such that
\[
\hat{u}-\hat{u}_\infty \in L^2(\hat{E}_1).
\] 
Above, $\hat{u}_\infty$ denotes the section  in $L^2_{loc}(\hat{E}_1)$ which is identically zero on $M_1(0)$ and coincides with the translation invariant  section $u_\infty$ on the infinite cylinder ${\bR}_+\times N$.   $u_\infty$ is uniquely determined  by $\hat{u}$ and thus we get  well defined map
\[
T_\infty: L^2_{ex}(\hat{E}_1)\ni \hat{u} \mapsto u_\infty \in L^2(E)
\]
called  {\em asymptotic limit (trace)} map.  $L^2_{ex} (\hat{E}_1)$ is naturally equipped  with a  norm
\[
\|\hat{u}\|^2 = \|\hat{u}-\hat{u}_\infty\|^2_{L^2(\hat{E}_1)}+\|u_\infty\|^2_{L^2(E)}.
\]
Clearly $L^2_{ex}(\hat{E}_1)$ with the above norm  is a Hilbert space and we have a short exact sequence 
\[
0 \hookrightarrow L^2(\hat{E}_1) \hookrightarrow L^2_{ex}(\hat{E}_1) \stackrel{T_\infty}{\ra}L^2(E)\ra 0.
\]
The map 
\[
L^2(E) \ni u_\infty \mapsto \hat{u}_\infty \in L^2_{ex}(\hat{E}_1)
\]
defines a splitting of this sequence.

Define the space of extended $L^2$-solutions of  $\hat{D}_{i, \infty}$ as 
\[
K_{i}:= \ker \hat{D}_{i, \infty} \cap L^2_{ex}(\hat{E}_i).
\]
The results of \cite{APS1}  and \cite{LM} show that these are finite dimensional spaces  and the spaces of asymptotic traces $L_i:=T_\infty(K_i)$ are subspaces in ${\cal H}:=\ker D\subset L^2(E)$. We have a difference map
\[
\Delta :K_1 \oplus K_2 \ra {\cal H}, \;\;\;\hat{u}_1\oplus \hat{u}_2 \mapsto T_\infty(\hat{u}_1)-T_\infty(\hat{u}_2).
\]
We denote its kernel by ${\cal K}_\infty$. It is a  finite dimensional subspace  of $L^2_{ex}(\hat{E}_1)\oplus L^2_{ex}(\hat{E}_2)$.

Finally we define the {\em splitting map}
\[
S_r: C^\infty({\cal E}_r) \ra L^2_{ex}(\hat{E}_1)\oplus L^2_{ex}(\hat{E}_2)
\]
by $\psi \mapsto S_r^1\psi \oplus S_r^2\psi$ where
\[
S_r^1 \psi =\psi \;\;{\rm on}\;\;  M_1(r)\subset M_1(\infty)
\]
and 
\[
(S^1_r\psi)(t,x) = \psi(r,x),\;\;\forall t\geq r,\; x\in N.
\]
$S_r^2$ is defined similarly. We  can now formulate the main result of this paper. When $\hat{B}_i\equiv 0$, i.e. the operators $\hat{D}_{i,\infty}$ are translation invariant along the neck, this result was  proved by  Cappell-Lee-Miller in \cite{CLM}   and is implicitly contained in \cite{Yos}. The super-symmetric situation was discussed in \cite{Mrow} in the special case of anti-selfduality operators.

\medskip

\noindent {\bf Main Theorem} {\em  Fix a positive  real number $\delta$ such that 
\[
\delta < \min (\gamma, \lambda)
\]
where  $\gamma$ is the smallest positive eigenvalue of $D$.  For every   function  $c: {\bR}_+\ra {\bR}_+$ such that
\[
c(r) =o(1/r),\;\; c(r)\geq Ce^{-\delta r} \;\;{\rm as}\;r\ra \infty
\]
 denote  by $\tilde{\cal K}_r(c)$  the subspace $L^2({\cal E}_r)$ of  spanned by the eigenvectors of ${\dir}_r$  corresponding  to eigenvalues  $|\lambda|\leq c(r)$. Then the  following hold.

\noindent(a) The splitting maps define an asymptotic map
\[
\tilde{\k}_r(c(r))\Lra^a K_1\oplus K_2.
\]
(b) The sequence
\[
0\ra \tilde{\k}_r(c(r))\stackrel{S_r}{\Lra^a}K_1\oplus K_2\stackrel{\Delta}{\ra}{\cal H}\ra {\cal H}\ra {\cal H}/(L_1+L_2)\ra 0
\]
is asymptotically exact. Furthermore, if all the operators involved are super-symmetric, the above  sequence is  super-symmetric as well.}

\medskip

In the course of the proof we will construct an asymptotic  inverse $\Psi_r$ to the splitting map which we call the {\em glueing map}.  This is an asymptotic map ${\k}_\infty\ra \tilde{\k}_r$ such  if $P_r$ denotes the orthogonal projection onto $\tilde{\k}_r$ and $P_\infty$ the orthogonal projection onto $\k_\infty$ then
\[
\|(P_\infty S_r)\circ (P_r \Psi_r) -{\bf id}_{\k_\infty}\|+\|(P_r \Psi_r)\circ (P_\infty S_r)-{\bf id}_{\tilde{\k}_r}\|=o(1) \;\;{\rm as}\;\; r\ra \infty.
\]

The rest of the paper is occupied with the proof and applications of the main theorem. It is  divided as follow. In Section 1 we list    some basic analytical  facts about elliptic equations on  manifolds with cylindrical ends.  We mention in particular the  {\bf Key Estimate}  which   adds a bit of compactness to the situation. Its proof is deferred to the Appendix.   Section 2 contains the proof of the Main Theorem itself. The strategy is  similar to the  approaches in \cite{CLM} and \cite{Yos}  but the details are  greatly simplified. Section 3 is devoted to a comparative study between the Main Theorem  and the  Mayer-Vietoris theorem for complexes of differential operators described in \cite{Andre}.  In the case of Dirac operators, we actually have a version of Poincar\'{e} lemma  which follows from the  existence results of \cite{BW}.    To construct the connecting homomorphism we follow closely its description in the DeRham case  contained in \cite{BT}.  This leads to a natural asymptotic connecting morphism. The  Main Theorem is equivalent with  an asymptotic Mayer-Vietoris sequence; see (\ref{eq: mv}).

In  Section 4 we present  two applications. The first one  is concerned with small eigenvalues  of selfadjoint  elliptic operators on manifolds containing long necks of the type considered by W.Chen  in \cite{Chen}. He  proved that if  ${\k}_\infty =0$ then the operators ${\dir}_r$  have no kernel  for $r \gg 0$  the norms of their inverses  are  $O(r)$.  We consider next  super-symmetric operators and we study what happens  if the space ${\k}_\infty$ is purely even.  We show that the  component  $C^\infty({\e}_r^+)\ra C^\infty({\e}_r^-)$ of $\dir_r$ admits  $L^2$-bounded right inverses of norms $O(r)$ as $r\ra \infty$. Such a situation is often encountered in  gauge theoretic gluing problems over {\em even} dimensional manifolds.

The second application, also suggested by problems in gauge theory, has to do  with families of operators. We  describe  ways to glue   the indices of some families of elliptic problems. 

We believe the asymptotic language  will find applications in other problems involving adiabatic deformations. It is not difficult to introduce the notion of asymptotic (co)chain complexes and asymptotic cohomology.  Many of  the basic  results in homological algebra have an asymptotic counterpart.

\bigskip

\noindent {\bf Acknowledgments} \hspace{.3cm} I am indebted to Ian Hambleton for asking the right questions which provided the stimulus for  writing this paper. Also, I want to thank  Ulrich Bunke and Matthias  Lesch who, back in 1994, pointed to me the fine resolution (\ref{eq: fine}).

%% file: spli1.tex
In this section we will  survey a few analytical facts which  are needed in the proof of the Main Theorem. The adequate functional background will be that of the Sobolev spaces $L^{k,p}$ consisting of distributions $k$-times differentiable with derivatives in $L^p$.

For any  $L^2_{loc}$ distribution $\hat{u}: t \mapsto u(t)$ on a cylinder $[0,L)\times N$ (where $L$ can be $\infty$) we denote by $\rho_t(\hat{u})$ the function $[0,L) \ra {\bR}_+$  defined by
\[
t\mapsto \rho_t(\hat{u}):=\left(\int_{C(t)}|u|^2 d\,vol\right)^{1/2},\;\;C(t)=[t,t+1]\times N.
\]
Additionally, define
\[
q:[0,L)\ra [0, \infty], \;\;t\mapsto q_{t,L}(\hat{u})= \sup_{t<s<L}\rho_{s}(\hat{u}).
\]
Note that if finite, $q_{t,L}$ is a decreasing function and thus  belongs to $L^\infty_{loc}(0,L)$.

Now let us observe that   the operator $J$ induces a symplectic structure on $L^2(E)$ defined by
\[
\omega(u,v):=\int_{N}(Ju, v) d\,vol
\]
The spectrum of $D$ is real and consists only of discrete eigenvalues with finite multiplicities.  Set
\[
{\cal H}_\mu:=\ker (\mu- D).
\]
and denote by $P_\mu$ the orthogonal projection onto ${\cal H}_\mu$. Since $\{J, D\}=0$ we deduce $J{\cal H}_\mu ={\cal H}_{-\mu}$.  The spectral gap of $D$ is the positive real number $\gamma=\gamma(D)$ defined as the smallest positive eigenvalue of $D$. Note that due to the spectral symmetry, $-\gamma(D)$ is also an eigenvalue of $D$. In particular, ${\cal H}={\cal H}_0$ is $J$ invariant  and   thus has an induced  symplectic structure. We have the following result (see  \cite{CLM}, \cite{Chen}, \cite{N}).

\begin{lemma} {\rm  The spaces  $L_i=T_\infty(K_i)$ of  asymptotic traces of  extended $L^2$ solutions are  lagrangian subspaces  of ${\cal H}$ i.e. $L_i^\perp=JL_i$.}
\end{lemma}

In the super-symmetric case we have  $E=E^+\oplus E^-$ and $G^*G ={\bf 1}_{E^+}$, $GG^*={\bf 1}_{E^-}$
\[
J=\left[
\begin{array}{cc}
0 & -G^* \\
G & 0
\end{array}
\right].
\]
Then  $J(E^\pm)=JE^\mp$ and 
\[
D= \left[
\begin{array}{cc}
{\bf D} & 0 \\
0 & J {\bf D}J^{-1}
\end{array}
\right].
\]
The space ${\cal H}$ is ${\bZ}_2$-graded
\[
{\cal H}={\cal H}^+\oplus {\cal H}^-
\]
and $G{\cal H}^+={\cal H}^-$. The asymptotic limit spaces $L_i$ now have decompositions
\[
L_i=L_i^+ \oplus L_i^-,\;\; L_i^\pm \subset {\cal H}^\pm
\]
and the lagrangian condition translates into
\be
(L_i^+)^\perp=G^*L_i^-, \;\;(L_i^-)^\perp =GL_i^+
\label{eq: lagr}
\ee
where $\perp$ denotes the orthogonal complements in ${\cal H}^\pm$.

Consider  a cylinder $[0,L) \times N$. Denote  by $\hat{E}$ the pullback of $E\ra N$ to this cylinder and by $\hat{D}$ the partial differential operator on $C^\infty(\hat{E})$  
\[
\hat{D}=\partial_t -D.
\]
For any   eigenvalue $\mu$ of $D$ and any smooth section $\hat{u}$ of $\hat{E}$  define a new section $\hat{u}_\mu$ by the condition
\[
\hat{u}_\mu\!\mid_{t\times N} =P_\mu u(t)\;\;\;u(t):=\hat{u}\!\mid_{t\times N}.
\]
Clearly $\hat{u}_\mu$ is a smooth section which we will regard  as a smooth  map
\[
\hat{u}_\mu= u_\mu(t): [0,L)\ra {\cal H}_\mu.
\]
Set
\[
u^\perp(t)=u(t)-u_0(t).
\]
\begin{proposition}{\bf (Key Estimate)}  {\rm There exists a constant $C>0$ depending (continuously) only  on the geometry of $N$  and the coefficients of $D$  with the following property. For any smooth  sections $\hat{u}$, $\hat{f}$ of $\hat{E}$ such that
\[
\hat{D}\hat{u}= \hat{f}
\]
and
\be
q_t(f) <\infty
\label{eq: density}
\ee
the following inequalities hold
\be
\|u_0(t)-u_0(t+n)\|_{L^2}\leq C\int_t^{t+n}q_{s,L}(\hat{f})ds, \;\;\;\forall  n\in {\bZ}\cap [0, L-t)
\label{eq: key1}
\ee
\be
\rho_{t+n}(\hat{u}^\perp)\leq C( e^{-\gamma n}\rho_t(\hat{u}^\perp) + e^{-\gamma n}\rho_{t+2n}(\hat{u}^\perp) + \frac{1}{\gamma^2}q_{t,L}(\hat{f})\,),\;\;\forall n\in {\bZ}\cap [0, (L-t)/2).
\label{eq: key2}
\ee
Above, $\gamma$ denotes the spectral gap of  $D$.}
\label{prop: key}
\end{proposition}

We have the following immediate consequence   whose   proof  is left to the reader.

 \begin{corollary}{\rm  Let $L=\infty$  in the} {\bf Key Estimate}.  {\rm If 
\[
\hat{D}\hat{u}= \hat{f}
\]
where both $\hat{u}$  and $\hat{f}$ are smooth and satisfy
\be 
\rho_t(\hat{u})\in L^\infty({\bR}_+),\;\;q_t(f)\in L^1({\bR}_+)
\label{eq: apriori}
\ee
then 
\[
\hat{u}\in L^2_{ex}(\hat{E})
\]
and}
\[
\|T_\infty\hat{u}-u_0(t)\|_{L^2}\leq C(\gamma^{-2}\int_t^\infty q_s(\hat{f}) ds+e^{-\gamma t}q_t(\hat{u})\,).
\]
\label{cor: key1}
\end{corollary}

Suppose $\hat{A}$ is a  smooth selfadjoint endomorphism of  $\hat{E}\ra {\bR}_+\times N$ such that  for some $ \lambda >0$ we have
\be
\sup\{|A(t,x)|\; ;\; x\in N\} =O(e^{-\lambda t}).
\label{eq: decay}
\ee
Set $\hat{A}_r:=\eta_r(t)\hat{A}$.

\begin{proposition}{\rm   Suppose  that we have a  sequence of  smooth sections $\hat{u}_r$ satisfying the following conditions.

\noindent (a) There exists $C >0$ such that $\rho_t(\hat{u}_r )< C$ for all  $t, r>0$.

\noindent (b)   The sections $\hat{u}_r$ and their derivatives are uniformly bounded on $C(0)$.

\noindent (c) There exists a sequence of smooth endomorphisms $B_r$ of $E$ such that
\[
m(r):=\sup\{|B_r(x)|\; ; \; x\in N\} = o(1/r)\;\;{\rm as}\; r\ra \infty
\]
 and $\hat{D}-\hat{A}_r\hat{u}_r  -B_r \hat{u}_r=0$ on the cylinder  $[0,r]\times N$.

\noindent (d)  $u_r(t) = u_r(r)$,  $\forall t\geq r\geq 0$.

Then    a subsequence of $\hat{u}_r$ converges {\em in the norm of}  $L^2_{ex}$ to a section $\hat{u}$  satisfying  $\hat{D}-\hat{A}\hat{u}=0$ on ${\bR}_+\times N$. Moreover, on a subsequence}
\be
u_r(r)_r \ra T_\infty \hat{u} \;\;{\rm in\;the \; norm \;of}\; L^2(E).
\label{eq: converg}
\ee
\label{prop: key2}
\end{proposition}

\noindent {\bf Proof}   In the sequel we will use the same symbol $C$ to denote positive constants independent of $t, r>0$. Set $\hat{f}_r=\hat{A}_r\hat{u}_r+B_r\hat{u}_r$.   Then
\be
\hat{D}\hat{u}_r= \hat{f}_r \;\;{\rm on}\; [0,r]\times N.
\label{eq: 1}
\ee
The conditions (\ref{eq: decay}), (a) and (b)  coupled with a standard bootstrap argument imply that there exists a constant  $C>0$ such that
\be
\sup\{|\hat{u}_r(t,x)|\; ;\; (t, x)\in [0,r-1]\times N\} \leq C,\;\;\forall r>0.
\label{eq: bound}
\ee
This implies that     a subsequence of $\hat{u}_r\!\mid_{[0,r]\times N}$ converges  weakly in $L^2_{loc}$ to a  section $\hat{u}$ defined over the entire cylinder.  Clearly $B_r \hat{u}_r \ra 0$ in $L^2_{loc}$ so that $\hat{u}$  is a weak solution  of
\[
\hat{D}\hat{u}-\hat{A}\hat{u}=0\;\;{\rm on}\; {\bR}_+\times N.
\]
We can now conclude   via elliptic estimates that  we can extract a subsequence which converges    is strongly  in $L^{k,2}_{loc}$. Moreover, according  to (\ref{eq: bound}) we deduce  $\hat{u} \in L^\infty$. If we now set $\hat{f}=\hat{A}u$  we deduce 
\[
\rho_t(f) \leq \|\hat{u}\|_\infty \rho_t(\hat{A}) = O(e^{-\lambda t}).
\]
Corollary \ref{cor: key1} implies $\hat{u}\in L^2_{ex}$ and
\be
\|u(t)-T_\infty\hat{u}\|_2 \leq C(e^{-\gamma t} +e^{-\lambda t}).
\label{eq: 2}
\ee 
The Key Estimate  for (\ref{eq: 1}), where 
\[
q_{t,r}(\hat{f}_r) \leq C (rm(r) +q_t(A_r) )\leq C(rm(r)+e^{-\lambda t}),\;\;r\geq t\geq 0.
\]
implies that  for all $0\leq t \leq r$ we have
\be
\|u_r(t)-u_r(r)\|_2\leq C (rm(r) +e^{-\lambda t}).
\label{eq: diff}
\ee
This proves  (\ref{eq: converg}) since $rm(r)=o(1)$.  To show that the convergence $\hat{u}_r\ra \hat{u}$  also takes place in the norm of $L^2_{ex}$ we only need  to establish  that on a subsequence
\[
\lim_{r\ra \infty}\int_0^\infty dt \int_N | (\,u_r(t)-u(t)\,)-(u_r(r)-T_\infty \hat{u})|^2 d\,vol  \ra 0.
\]
We extract the subsequence using the following argument.  For every  $n>0$    pick  $r=r_n>n$ such that the following inequalities hold.
\be
\int_0^ndt \int_N |u_r(t)-u(t)|^2 d\, vol \leq  \frac{1}{n^2}
\label{eq: 3}
\ee
\be
\int_{N}|u_{r_n}(n)-u_{r_n}(r_n)|^2 d\, vol < \frac{1}{n^2}
\label{eq: 4}
\ee
\be
\int_n^\infty dt \int_N|u(t)-T_\infty\hat{u}|^2 d\, vol \leq \frac{1}{n^2}
\label{eq: 5}
\ee
The choice (\ref{eq: 3}) is possible  because  the  sequence $\hat{u}_r$ converges  to $\hat{u}$ in the norm $L^2([0, n]\times N)$. The choice (\ref{eq: 4}) is possible because  $rm(r)=o(1)$ and (\ref{eq: diff}). Finally,  the choice (\ref{eq: 5}) is possible because of (\ref{eq: 2}). The subsequence $\hat{u}_{r_n}$ chosen as above converges  to $\hat{u}$ in the norm of $L^2_{ex}$.   $\Box$

%% file: spli2.tex
To show that $\lim_{r\ra \infty}\delta(S_r(\tilde{\k}_r), {\k}_\infty)=0$ we will use the following elementary result which follows immediately from Lemma \ref{lemma: asym}. 

\begin{lemma}{\rm Suppose $U$ is a finite dimensional subspace in a Hilbert space and $U_r$ is a sequence of finite dimensional  subspaces such that
\be
\lim_{r\ra \infty}\hat{\delta}(U_r, U)=0.
\label{eq: gap}
\ee
and 
\be
{\rm liminf}\, \dim U_r \geq  \dim U.
\label{eq: dim}
\ee
Then}
\[
\lim_{r\ra \infty} \delta(U_r, U)=0.
\]
\label{lemma: gap}
\end{lemma}

We will show that the two assumptions in the lemma are satisfied   if $U_r=S_r\tilde{\k}_r$ and $U={\k}_\infty$.  The proof of  the Main Theorem is thus  divided in two  steps.

\noindent {\bf Step 1}
\[
\lim_{r\ra \infty}\hat{\delta}(S_r(\tilde{\k}_r), {\k}\infty)=0.
\]
We argue by contradiction.  Thus  we assume there exists a sequence $\psi_r \in \tilde{\k}_r$  such  that
\be
\|S_r\psi_r\|_{ex}=O(1) \;\;{\rm as}\; r\ra \infty
\label{eq: 21}
\ee
and there exists $d_0>0$ such that
\be
{\rm dist}\, (S_r\psi_r , {\k}_\infty) >d_0,\;\;\forall r>0.
\label{eq: 22}
\ee
Set $\psi_r^i:=S_r^i\psi_r$, $i=1,2$.  We study only the behavior of $\psi_r^1$.  The sequence $\psi_r^2$ behaves similarly.   The condition (\ref{eq: 22}) shows   there exists a constant $c>0$ such that 
\[
\|\psi_r^1\|_{ex} \geq c\;\;\forall r>0.
\]
Thus we can normalize $\psi_r^1$ so that $\|\psi_r^1\|_{ex}=1$ and  (\ref{eq: 22}) continues to hold (with an eventually smaller $d_0>0$).

Note  first that using standard elliptic estimates  and (\ref{eq: 21})  we deduce  that $\psi_r^2$ and its derivatives are uniformly bounded on  $M_1(0)$. Thus   a subsequence of $\psi_r^1$ converges to a solution of $\hat{D}_{1,\infty}\hat{u}=0$ on $M_1(0)$. Using  Proposition \ref{prop: key2} we deduce  that  a further subsequence of the restriction of $\psi_r^1$ to ${\bR}_+\times N$ converges  in the norm of $L^2_{ex}$ to a solution of $\hat{D}_{1,\infty}\hat{\psi}=0$ on this  semi-infinite cylinder.  Clearly  we have produced a    section $\psi_\infty^1\in \ker \hat{D}_{1,\infty}\cap L^2_{ex}$ of norm $1$.    We proceed similarly with $\psi_r^2$. We now have  a pair
\[
\Psi:=\psi_\infty^1\oplus \psi_\infty^2 \in K_1\oplus K_2
\]
of norm $2$ which according to (\ref{eq: converg}) in Proposition \ref{prop: key2} have the same asymptotic limit.  Thus $\Psi \in {\k}_\infty$.  However, this contradicts  (\ref{eq: key2}).  Step 1 is completed.

\medskip

\noindent {\bf Step 2}     We will prove that 
\[
\dim {\k}_\infty \leq \dim S_r \tilde{\k}_r\;\;\forall r\gg 0.
\]
We will rely on the following auxiliary result.

\begin{lemma}{\rm Suppose $u\in L^{1,2}({\e}_r)$ is such that
\[
\|\dir_ru\|_2 < (1-\ve)c(r) \|u\|_2.
\]
Then  ${\rm dist}\,(u, \tilde{\k}_r(c(r))) < (1-\ve)\|u\|_2$.}
\label{lemma: dist}
\end{lemma}

\noindent {\bf Proof of the lemma}\hspace{.3cm}  Using the orthogonal decomposition
\[
L^2({\e}_r)=\tilde{\k}_rc(r) \oplus \tilde{\k}_r(c(r))^\perp
\]
we can write
\[
u=v +v^\perp.
\]
Then  ${\rm dist}\,(u, \tilde{\k}_r)=\|v^\perp\|_2$.  On the other hand
\[
(1-\ve)^2c(r)^2\|u\|^2 > \|{\dir}_ru\|^2  \geq \|{\dir}_rv^\perp\|^2  \geq  \Lambda^2 \|v^\perp\|^2
\]
where $\Lambda^2 > c(r)^2$. The lemma is proved.   $\Box$

\bigskip

To conclude the proof of Step 2 we will construct  for  $r\gg 0$ a space  $V_r\subset  L^2({\e}_r)$ isomorphic to ${\k}_\infty$ such that
\be
\hat{\delta}(V_r, \tilde{\k}_r) <1.
\label{eq: 23}
\ee
According to Lemma 2.3, Chap. IV, {\S}2 in \cite{Kato} this means that the orthogonal projection  onto  $\tilde{\k}_r$ induces an injection  $V_r \ra \tilde{\k}_r$ so that
\[
\dim {\k}_\infty =\dim V_r \leq \dim \tilde{\k}_r,\;\;\forall r\gg 0.
\]
The condition (\ref{eq: 23}) is satisfied provided  ${\rm dist}\,(v, \tilde{\k}_r) < v$, for all $v\in V_r\setminus\{0\}$. According to Lemma \ref{lemma: dist} is suffices to construct a subspace  $V_r\subset L^{1,2}({\e}_r)$ isomorphic to  ${\k}_\infty$ such that
\be 
\sup_{v\in V_r\setminus\{0\}}\frac{\|{\dir}_rv\|_2}{\|v\|_2} < c(r).
\label{eq: 24}
\ee
Such a subspace is obtained via a simple glueing construction.

We construct a glueing map 
\[
\Psi_r: {\k}_\infty \ra L^{1,2}({\e}_r),\;\;\hat{u}_1\oplus \hat{u}_2 \mapsto \Psi_r
\]
 uniquely determined by the following conditions.  Let $u_\infty$ denote the common asymptotic limit of  $\hat{u}_i$.  Now set
\[
\hat{v}_1= \eta_r (t)\hat{u}_1 +(1-\eta_r(t)) u_\infty.
\]
Define $\hat{v}_2$ similarly. Clearly on the overlap
\[
\imath_r:[r+1, r+2]\times N \ra [-r-2,-r-1] \times N
\]
we have
\[
\hat{v}_1=\hat{v}_2=u_\infty
\]
 so we  can glue these two sections on the overlap to produce  a  smooth section $\Psi_r\in C^\infty ({\e}_r)$. Clearly the map $\Psi_r$ is linear. Set $V_r:=\Psi_r({\k}_r)$.

Note that $\Psi_r$ is injective because if $\Psi_r(\hat{u}_1,\hat{u}_2)\equiv 0$ then both $\hat{u}_i$ must vanish on $M_i(0)$ and by unique continuation they must vanish everywhere. We claim $V_r$ satisfies (\ref{eq: 24}).

Clearly  ${\dir}_r \Psi_r\equiv 0$ on $M_1(r-1), M_2(-r+1)\subset M(r)$ so we only need an   estimate of $\hat{D}_{1,r}\hat{v}_1$ on the cylinder  $[r-1,r+2]\times N$ and a similar one for $\hat{D}_{2,r}\hat{v}_2$. 

On this cylinder we have $\hat{D}_{1,r}=J(\partial_t -D-\eta_r\hat{A}_1)$ so that we have
\[
-J\hat{D}_{1,r}\hat{v}_1=(\,\hat{D}-\eta_r\hat{A}_1\, )(\, \eta_r\hat{u}_1+(1-\eta_r)u_\infty\, )
\]
\[
=  (\, \hat{D}-\hat{A}_1\, )  (\, \eta_r\hat{u_1}+(1-\eta_r)u_\infty \,)  
\]
\[
+  ( 1-\eta_r)\hat{A}_1(\eta_r\hat{u_1}+(1-\eta_r)u_\infty ) .
\]
We examine the two  terms separately. The first one can be rewritten as
\[
(\hat{D}-\hat{A}_1 )  (\eta_r\hat{u_1}+(1-\eta_r) u_\infty )  = [ \hat{D}-\hat{A}_1, \eta_r ]\hat{u}_1 +[ \hat{D}-\hat{A}_1, (1-\eta_r) ] u_\infty  
\]
\[
+\eta_r (\hat{D}-\hat{A}_1)\hat{u}_1+ (1-\eta_r)(\hat{D}-\hat{A}_1)u_\infty
\]
($\hat{c}$=Clifford multiplication on ${\bR}_+\times N$)
\[
= \hat{c}(d\eta_r)(u_1(t)-u_\infty)  - (1-\eta_r)\hat{A}_1u_\infty.
\]
Thus 
\[
-J\hat{D}_{1,r}\hat{v}_1= \hat{c}(d\eta_r)(u_1(t)-u_\infty) +\eta_r( 1-\eta_r)\hat{A}_1(\hat{u}_1-u_\infty).
\]
We can now use  Corollary \ref{cor: key1} for the equation
\[
\hat{D}\hat{u}_1 =\hat{f}=\hat{A}_1\hat{u_1}\;\;{\rm on}\;\;{\bR}_+\times N.
\]
 The   decay rate of $A_1(t)$  shows that $\rho_t(\hat{f}) =O(e^{-\lambda t}\rho_t(\hat{u}_1)\,)$. Hence
\be
\|{\dir}_r\hat{v}_1\|_{2,[r-1,r+2]\times N} \leq C (e^{-\lambda r} + e^{-\gamma r}) q_t(\hat{v}_1).
\label{eq: 26}
\ee
(\ref{eq: 26}) implies (\ref{eq: 24})  since $q_t(v_1) \leq C \|\Psi_r\|$  (for $r\gg 0$) and 
\[
c(r)> C(e^{-\lambda r}+ e^{-\gamma r}), \;\;\forall r\gg 0.
\]
The Main Theorem is proved  $\Box$

\medskip

\begin{remark}{\rm The exponetial decay condition (\ref{eq: decay1}) on the coefficients of $\hat{D}_{i,\infty}$ can be replaced by a  milder one
\[
\sup \{ |\hat{A}_i(\hat{x})|\, ;\; \hat{x}\in [t,t+1]\times N\} \leq C t^{-p}
\]
where  $p>2$.  The statement of the Main Theorem  changes in an obvious way to  take this decay into account.}
\end{remark}

\begin{remark}{\rm We can rewrite the  conclusion of the Main Theorem as a short asymptotically exact sequence
\[
0\ra \tilde{\k}_r\stackrel{S_r}{\Lra^a}K_1\oplus K_2\stackrel{\Delta}{\ra} L_1\oplus L_2 \ra 0.
\]
The gluing map  $\Psi_r$ is an asymptotic splitting of this sequence, in the sense described in the introduction.}
\end{remark}

\begin{remark}{\rm The Main Theorem extends easily to families of operators.   Suppose $X$ is a compact  CW-complex and   all the constructions in  the introduction depend continuously on the parameter $x\in X$ such that   the  spectral gaps of the boundary operators $D_x$ are bounded  from below
\[
\gamma_0:=\inf_{x\in X} \gamma(D_x) >0.
\]
Then $h(x):=\dim {\cal H}_x$ is independent of $x$. We denote this common dimension by $h$. Assume also  that  the  functions
\[
\kappa_i:X \ra {\bZ} ,\;\;\kappa_i(x):= \dim K_i(x)\;\;(i=1,2)
\]
\[
\ell: X\ra {\bZ},\;\;\ell(x)=\dim L_1(x)\cap L_2(x)
\]
are constant, $k_i(x)\equiv \kappa_i$, $\ell(x)\equiv \ell$. One then can show that  $K_i(x)$ and $L_1(x)\cap L_2(x)$ depend continuously upon $x$ in the gap topology. Thus they can be viewed as  continuous maps  in grassmannians of finite dimensional subspaces  in    Hilbert spaces and as such they define  vector bundles over $X$.   The Main Theorem   for families  states that for $r\gg 0$  the spaces $\tilde{\k}_{r,x}(c)$ form a vector bundle over $X$ and we   have  and exact sequence of vector bundles
\[
0\ra\tilde{\k}_r \stackrel{\Gamma_r}{\ra} K_1\oplus K_2 \stackrel{\Delta}{\ra}{\cal H}\ra {\cal H}/(L_1 +L_2) \ra 0.
\]
Since $L_1$, $L_2$ are lagrangian then
\[
{\cal H}/(L_1+L_2) \cong (L_1+L_2)^\perp =L_1^\perp\cap L_2^\perp=J(L_1\cap L_2).
\]
A similar statement is true in the super-symmetric case. In \cite{N2} we    describe a general gluing formula  for index of families when the   functions $h(x)$, $\ell(x)$ and $\kappa_i(x)$ are not necessarily constant.

In practice, one often encounters a {\em fibered version} of the problem. Suppose $U$, $V$ and $Z$ are  compact CW complexes and  $\zeta_U:U\ra Z$, $\zeta_V:V\ra Z$ are continuous maps.  We think of $U$ as a parameter space for a continuous family of operators 
\[
\hat{D}_{1,\infty}(u)=J(\partial_t -D_u-A_{1,u})
\]
 on $\hat{E}_1$ while $V$ is a parameter space for a continuous family of operators 
\[
\hat{D}_{2,\infty}(v)=J(\partial_t -D_v-A_{2,v})
\]
on $\hat{E}_2$.  $Z$ parameterizes a continuous family of elliptic operators $D_z$ on $E$ such that  $D_z=D_u=z=D_v$ if $z=\zeta(u)=\zeta(v)$.    Form the fiber product
\[
X=\{ (u,v)\in U\times V\; ;\; \zeta(u)=\zeta(v)\} =\zeta^{-1}(\Delta_Z)
\]
where $\Delta_Z$ is the diagonal $\{(z,z)\in Z\times Z\; ; \; z\in Z\}$.  For each $x=(u,v)\in X$ we  set $\hat{D}_{1,\infty}(x):=\hat{D}_{2,\infty}(x)$. It is precisely this fibered context  used  in \cite{Mrow} to establish the  super-symmetric version of the Main Theorem for families of anti-selfduality operators.}
\label{rem: families}
\end{remark}

%% file: spli3.tex
 The Main Theorem has a Mayer-Vietoris flavor.   One can formulate a genuine Mayer-Vietoris theorem as  follows. Consider the super-symmetric Dirac operator ${\dir}_r$ on ${\e}_r$. For simplicity we drop the  subscript $r$ and we assume the operators $\hat{D}_{i,\infty}$ are actually translation invariant. The results below  do hold for exponentially decaying perturbations as well.  The operator ${\dir}={\dir}_r$ has a  block decomposition
\[
\dir=\left[
\begin{array}{cc}
0 & \cdir^* \\
\cdir & 0 
\end{array}
\right]
\]
Consider now the  sheaf ${\cal S}$ on $M=M(r)$     defined by
\[
{\cal S}(u):= \{u\in C^\infty({\e}^+\!\mid_U)\; ;\; \cdir u=0\}
\]
for any open $U\subset M$.  Denote by $\Gamma^\pm$ the sheaves of  locally defined smooth sections of ${\e}^\pm$. Then the sequence
\be
0\ra {\cal S}\ra  \Gamma^+ \stackrel{{\cdir}}{\ra}\Gamma^- \ra 0
\label{eq: fine}
\ee
 is a fine resolution of ${\cal S}$.   Indeed, the exactness at $\Gamma^+$ is tautological. The exactness at $\Gamma^-$ is equivalent  to local existence of ${\cdir}u=f$, $\forall f\in \Gamma^-$. This can be proved easily using the existence and regularity results in \cite{BW}.  Thus  the Cech cohomology of the sheaf $M$ on $M$ is given by the complex
\[
(\Gamma, D)\;\;\;0\ra C^\infty({\e}^+)\stackrel{\cdir}{\ra}C^\infty({\e}^-)\ra 0
\]
so that
\be
H^0(M, {\cal S})\cong \ker \cdir,\;\; H^1(M, {\cal S})\cong {\rm coker}\, \cdir,\;\;H^k(M, {\cal S})=0,\;\;\forall k\geq 2.
\label{eq: coho}
\ee
In particular,
\[
\chi(M, {\cal S})= {\rm ind}\,(\cdir).
\]

\begin{remark}{\rm  The same techniques in \cite{N4} can be used to show that the element in $K$-theory defined by the symbol of $\cdir$ is completely determined by  the sheaf ${\cal S}$. This implies the index of $\cdir$ is completely determined by this sheaf. The equality (\ref{eq: coho})  explains how.}
\end{remark}

Denote by $U_i$ an  open neighborhood of $M_i(r)\subset M(r)$, slightly larger that $M_i(r)$. As in \cite{Andre} or \cite{BT} one can show that  we have a short  Mayer-Vietoris exact sequence of  co-chain complexes
\[
0\ra  (\Gamma(M), \cdir) \ra (\Gamma(U_1), \cdir)\oplus (\Gamma(U_2),\cdir) \stackrel{\Delta}{\ra}\ra \Gamma(U_2\cap U_2)\ra 0
\]
from which we get a long Mayer-Vietoris sequence
\[
0\ra \ker \cdir \ra {\cal S}(U_1)\oplus  {\cal S}(U_2)\ra {\cal S}(U_1\cap {\cal S}(U_2)\stackrel{\partial}{\ra} 
\]
\[
\stackrel{\partial}{\ra}  \ker \cdir^* \ra H^1(U_1, {\cal S})\oplus H^1(U_2, {\cal S})\ra H^1(U_1\cap U_2, {\cal S})\ra 0.
\]
Imitatting \cite{BT}, we can provide  quite an explicit description of  the connecting morphism. Choose a partition of unity subordinated to the cover $\{U_1, U_2\}$.  Denote it by    $\{\omega_i \in C_0^\infty(U_i)\}$. Let $u\in {\cal S}(U_1\cap U_2)$. The section $\omega_2 u$, extended by zero outside $U_2$ will be regarded as a smooth section over $U_1$ which we denote by $u_1$. Similarly, we can regard $-\omega_1 u$ as a smooth section $u_2$ over $U_2$. Note that
\[
\cdir (u_1-u_2)= \cdir u =0\;\;{\rm on}\;\; U.
\]
Thus the two sections $\cdir u_i$ match-up to define a section $\partial u \in C^\infty({\e}^-)$.   We can give $\partial u$ a more suggestive description using the  equalities
\[
\cdir u_1= [\cdir , \omega_2] u+\omega_2\cdir u=c(d \omega_2)u
\]
and similarly,
\[
\cdir u_2=-c(d\omega_1) u
\]
where $c(?)$ denotes the Clifford multiplication by a $1$-form.

This suggests the existence of an asymptotic map
\[
{\cal H}^+\Lra^a \ker \cdir_r^*
\]
defined by
\[
\partial:\;\;u\mapsto c(d\omega_2) u.
\]
To evaluate  how far is it from being in the kernel of $\cdir_r^*$ we use Proposition 3.45 in \cite{BGV} and we deduce
\[
\cdir^*c(d\omega_2)u =-c(d\omega)\cdir u  -2\nabla_{{\bf grad}\,\omega_2}u + (\Delta \omega_2)u = -2\nabla_{{\bf grad}\omega_2}u + (\Delta \omega_2)u
\]
Now choose $\omega_2$ so that it depends only on the longitudinal coordinate  $t$ along the neck.  In this case we deduce
\[
\cdir^*(c(d\omega_2) =-\ddot{\omega}_2 u
\]
where we use dots to denote $t$-derivatives.  The  quantity
\[
\frac{\|\cdir^*(c(d\omega_2)u\|^2}{\|c(d\omega_2)u\|^2} =\frac{|\ddot{\omega}_2|^2}{|\dot{\omega}_2|^2}
\]
offers an indication about the distance between $\partial u:=c(d\omega_2)u$ and $\tilde{\k}_r^-$.   It suggest  choosing $\omega_2$ such that $\dot{\omega}_2$ is close the first eigenvalue of the Dirichlet problem of the one dimensional Laplacian on an interval of length $\sim r$. In such an instance,  we see that the above quantity is $O(r^{-2})$.  For example ($L = L(r)= \kappa r$, $\kappa \in (0,1)$)
\[
\dot{\omega}_2=\frac{\pi}{2L}\alpha (t)\sin (\frac{\pi t}{L})
\]
where $\alpha(t)$ is a nonnegative, smooth cut-off function supported in $[0,L]$ such that $\alpha(t)\equiv 1$ on $[1,L-1]$ and there exists a constant $C>0$ independent of $L$ such that $\frac{d^k\alpha}{dt^k}(t)|\leq 2$, $\forall t$, $k=1,2,3$. Now set
\[
\omega_2(t) = \frac{\pi}{2L}\int_0^t\alpha(s)\sin(\frac{\pi s}{L}) ds.
\]
We denote  by $\partial_r u$  the orthogonal projection onto $\tilde{\k}_r^-$ of
 \[
c(d\omega_{2}(t+ L/2) )u \in C^\infty({\e}_r^-).
\]

\begin{definition}{\rm The sequence of asymptotic maps 
\[
U_r\stackrel{f_r}{\Lra^a}V_r\stackrel{g_r}{\Lra^a}W_r
\]
is said to be {\em weakly asymptotically exact} if for $r\gg 0$ there exists an  exact sequence
\[
U_r\stackrel{f'_r}{\Lra^a}V_r\stackrel{g'_r}{\Lra^a}W_r
\]
such that}
\[
 \delta({\rm Graph}(f_r), {\rm Graph}(f'_r))+  \delta({\rm Graph}(g_r), {\rm Graph}(g'_r)) =o(1)\;\;{\rm as}\;\;r\ra \infty.
\]
\end{definition}

Note that Lemma \ref{lemma: asym} (see also Remark \ref{rem: wae}) shows that any asymptotically exact sequence is also weakly asymptotically exact.

One can show that we have a weakly  asymptotically  exact sequence
\[
0\ra \tilde{\k}^+_r \stackrel{S^+_r}{\Lra^a} K^+_1\oplus K_2^+ \stackrel{\Delta}{\ra} L^+_1+L_2^+ \stackrel{\partial_r}{\Lra^a}
\]
\be
\Lra^a \tilde{\k}_r^- \stackrel{S_r^-}{\Lra^a}K_1^-\oplus K_2 \stackrel{\Delta^-}{\ra} L_1^-\oplus L^-_2 \ra 0.
\label{eq: mv}
\ee
Moreover, $\ker S_r^- = 0$, $\forall r>0$. This sequence is only weakly asymptotically exact because the range of $\partial_r$ is never trivial  but nevertheless, 
\be
\|\partial_r\|^2= o(1). 
\label{eq: small}
\ee
The  above estimate  follows from the fact that $c(r) =o(r^{-1})$, $L\sim  r$ and
\[
\|\dir_r\phi_r -\frac{\pi}{L} \phi_r\|^2_{L^2(M(r))} =o(1) \| \phi_r\|^2_{L^2(M(r))}
\]
where 
\[
\phi_r= \ddot{\omega}_2(t+L/2)u \oplus \dot{\omega}_2(t+L/2)Gu \in C^\infty({\e}_r^+)\oplus C^\infty({\e}_r^-).
\]
We leave the details to the reader.

Formula (\ref{eq: mv}) predicts
\[
{\rm ind}\,(\cdir)= ( \dim K_1^+ -\dim K_1^- ) + ( \dim K_2^+ - \dim K_2^- )
\]
\be 
 - \dim ( L_1^+ + L_2^+ ) + \dim ( L_1^- + L_2^-) .
\label{eq: index}
\ee
This  equality can be alternatively established as follows.

Set $h=\dim {\cal H}^\pm$ and
\[
U_i^+ = \ker \hat{\bf D}_i \cap L^2(M_i(\infty ; \hat{E}^+_i)
\]
\[
U_i^-=\ker \hat{\bf D}_i^* \cap L^2(M_i(\infty); \hat{E}_i^-), \;\;i=1,2.
\]
We have short exact sequences
\[
0\ra U_i^\pm \ra K_i^\pm \ra L_1^\pm \ra 0.
\]
The Atiyah-Patodi-Singer index theorem  (Thm. 3.10 in \cite{APS1}) coupled with the Atiyah-Singer index theorem on closed manifolds  implies  immediately
\[
{\rm ind}\, \cdir = (\dim U_1 ^+ - \dim K_1^-)+(\dim U_2^+ -\dim K_2^-) + h
\]
\be
=  ( \dim K_1^+ -\dim K_1^- ) + ( \dim K_2^+ - \dim K_2^- )  +h -(\dim L_1^+ +\dim L_2^+).
\label{eq: index1}
\ee
Now observe that
\[
\dim L_1^+ +\dim L_2^+ = \dim (L_1^+ + L_2^+) +\dim (L_1^+\cap L_2^+)
\]
\[
=  \dim (L_1^+ + L_2^+) + h -\dim (\, (L_1^+)^\perp + (L_2^+)^\perp \,)
\]
(the above $\perp$ denotes the orthogonal complement in ${\cal H}^+$)
\[
\stackrel{(\ref{eq: lagr})}{=}  \dim (L_1^+ + L_2^+) + h - \dim (G^*(L_1^-+L_2^-)
\]
\[
= h +  \dim (L_1^+ + L_2^+) -\dim (L_1^- +L_2^-).
\]
Using this last equality in (\ref{eq: index1}) we obtain (\ref{eq: index}).

%% file: spli4.tex
As promised, we will include some simple applications of the Main Theorem.  

Suppose for example that ${\k}_\infty =0$.  This is possible  if and only if
\[
L_1\cap L_2 =0
\]
and  $\ker(T_\infty: K_i \ra L_i)=0$, $i=1,2$.  These kernels consist of the $L^2$-solutions of $\hat{D}_{i,\infty}$.     This shows that the operators ${\dir}_r$ cannot have eigenvalues $\lambda_r$ such that $|\lambda_r|=o(1/r)$ as $r\ra \infty$. We have thus established the following result (proved for the first time in \cite{Chen}).

\begin{corollary} {\rm  Suppose that
\[
L_1\cap L_2 =\{0\}\;\;{\rm and}\;\; \ker \hat{D}_{i,\infty} \cap L^2(\hat{E}_i)= \{0\}
\]
Then  for $r\gg 0$ the operator $\dir_r$ has  a bounded inverse
\[
\dir_r^{-1}:L^2({\e}_r)\ra L^2({\e}_r)
\]
and}
\[
\|\dir_r^{-1}\| =O(r)\;\;{\rm as}\; r\ra \infty.
\]
\end{corollary}

Suppose now the entire situation is super-symmetric.  Thus, we have  decompositions
\[
K_i=K_i^+\oplus K_i^-,\;\;{\k}_\infty={\k}_\infty^+\oplus {\k}_\infty^-.
\]
In \cite{Mrow},  $\tilde{\k}_r^-$ was called the {\em obstruction space}. We assume 
\be
K_i^-=\{0\},\;\;i=1,2.
\label{eq: 31}
\ee
This implies $L_i^-=\{0\}$ and   ${\k}_\infty^-=\{0\}$.  The equality (\ref{eq: lagr})   shows that $L_1^+=L_2^+={\cal H}^+$. We  deduce
\be
\tilde{\k}_r^-=\{0\},\;\;\forall r\gg 0
\label{eq: 32}
\ee
while the  even part $\tilde{\k}_r$ fits in an exact sequence
\[
0\ra \tilde{\k}_t \stackrel{\Gamma_r^+}{\ra}K_1\oplus K_2 \stackrel{\Delta^+}{\ra}{\cal H}^+\ra 0.
\]
The bundle ${\e}_r$ has a decomposition
\[
{\e}_r={\e}_r^+\oplus {\e}_r^-
\]
with respect to which ${\dir}_r$ has the super-symmetric block decomposition
\[
{\dir}_r=\left[
\begin{array}{cc}
0 & \not\!\!\dir_r^* \\
\not\!\!\dir_r & 0 
\end{array}
\right]
\]
where $\cdir_r:C^\infty({\e}_r^+)\ra C^\infty({\e}_r^-)$. The equality (\ref{eq: 32}) implies that $\cdir_r$ is onto since 
\[
\tilde{\k}_r=\ker\cdir_r^* \cong {\rm coker}\,\cdir_r.
\]
Thus $\cdir_r \cdir_r^*$ is one-to-one and onto and admits a bounded inverse $L^2({\e}_r^-)\ra L^2({\e}_r^-)$.  We claim that
\be
\|(\cdir_r\cdir_r^*)^{-1}\|=O(r^2) 
\label{eq: 33}
\ee     
To prove this claim we argue by contradiction.

Because $\cdir_r\cdir_r^*$ is self-adjoint, positive  and has  compact resolvent, the norm of its  inverse  is $m(r)^{-1}$ where
\[
m(r) = \inf\{ \lan \, \cdir_r\cdir_r^*u\, , \, u\, \ran ;\;\;\|u\|=1\}.
\]
 Suppose that for every $r\gg 0$ we can find $\phi^r\in L^{2}({\e}_r^-)$ such that
\[
\|\phi_r\|=1 
\]
and
\[
m(r)=\|\cdir_r^*\phi_r\|^2=(\cdir\cdir_r^*\phi_r, \phi_r)=o(1/r^2) \;\;{\rm as}\; \ra \infty.
\]
Now pick $c(r)  >  \exp(-\delta(r))$ such that
\[
c(r)=o(1/r),\;\; m(r)=o(c(r)^2)\;\;{\rm as}\; \;r\ra \infty.
\]
The above  $\delta$ is  the same exponent as in  the Main Theorem.

Now apply $\dir_r$ to the vector $u_r:=0\oplus \phi_r\in L^2({\e}_r^+\oplus{\e}_r^-)$. We deduce
\[
\|u_r\|=1
\]
and 
\[
\frac{\|{\dir}_r u_r\|}{\|u_r\|} =\sqrt{m(r)}.
\]
Thus, according to Lemma \ref{lemma: dist} we can conclude
\[
{\rm dist} \,(u_r, \tilde{\k}_r(c(r)))\leq  \frac{\sqrt{m(r)}}{c(r)}=o(1).
\]
On the other hand, $u_r$ is purely odd which implies ${\k}^-_r(c(r))\neq 0$ for all $r\gg 0$. This contradicts (\ref{eq: 32}) and  thus proves (\ref{eq: 33}).  This estimate also shows that $\cdir_r$ has a right inverse $R_r: L^2({\e}_r^-)\ra L^2({\e}_r^+)$ of norm $O(r)$.  We can now state our next result.

\begin{proposition}{\rm Suppose  the condition (\ref{eq: 31}) is satisfied. Then  for $r\gg 0$ the operator  $\cdir_r$  is onto and admits a   bounded right inverse of norm $O(r)$. Moreover}
\be
\ker \cdir_r =\tilde{\k}_r^+\;\;{\rm for}\;\;r\gg 0.
\label{eq: kernel}
\ee
\label{prop: inverse}
\end{proposition}

\noindent {\bf Proof}\hspace{.3cm} The only thing left to prove is the equality  (\ref{eq: kernel}) which follows immediately from the fact that the index of $\cdir_r$ is independent of $r$  and 
\[
\dim \tilde{\k}_r^+ -\dim\tilde{\k}_r^- ={\rm ind}\, \cdir = \dim \ker \cdir_r-\dim \ker \cdir_r^*.\;\;\Box
\]

\begin{remark}{\rm Results of this type are needed  in gauge theoretical gluing problems over ({\em even dimensional}) smooth manifolds. In such problems, the condition  (\ref{eq: 31})   appears in the following disguise.  

The operators $\hat{D}_{i,\infty}$ (defined in the introduction) have a super-symmetric decompositions
\[
\hat{D}_{i,\infty}=\left[
\begin{array}{cc}
0 & \hat{\bf D}^*_{i,\infty} \\
\hat{\bf D}_{i,\infty}
\end{array}
\right].
\]
For simplicity, we will omit the subscripts $i$, $\infty$. According to \cite{LM},  the operator $\hat{\bf D}$ induces  a Fredholm operator
\[
L^{1,2}_\delta(\hat{E}^+)\ra L_\delta^2(\hat{E}^-)
\]
where  $\delta$  is a small positive number      and the $L^{k,2}_\delta$ norm  is the $L^{k,2}$ norm with respect to  a measure on $M_i(\infty)$ which along the neck has the form $e^{\delta |t|} dt\wedge d\,vol_N$. The condition (\ref{eq: 31}) signifies that  $\hat{\bf D}$, in this functional set-up, is onto.}
\end{remark}

Suppose now that in Proposition \ref{prop: inverse} we have a  family of operators, each satisfying  (\ref{eq: 31}) and subject to the restrictions  listed in Remark \ref{rem: families}.  We deduce immediately the following  consequence.

\begin{corollary}{\rm   Under the above assumptions, there exists an exact sequence of vector bundles
\[
0\ra\ker \cdir_r \stackrel{\Gamma_r}{\ra} K_1\oplus K_2 \stackrel{\Delta^+}{\ra} {\cal H}^+\ra 0.
\]
In particular, by passing to determinant line bundles we deduce an isomorphism of line bundles over $X$
\[
\det ({\rm ind}\,(\cdir_r) \cong \det K_1 \otimes \det K_2 \otimes (\det {\cal H}^+)^*.
\]
where in the right hand side ${\rm ind}\,(\cdir_r)$ is viewed as an element in an appropriate $K$-theory of the parameter space $X$}
\label{cor: det}
\end{corollary}

\begin{remark}{\rm (a) The terms $\det K_i$  are also    determinant line bundles associated to  the indices of   the families of  Atiyah-Patodi-Singer problems  determined by $\hat{\bf D}_i$, $i=1,2$.    

(b) Corollary \ref{cor: det} is also useful in   orientability issues   involving various moduli spaces  arising in gauge theory.}
\end{remark}

%% file: splia.tex
The Key Estimate   is a consequence of the following elementary result.

\begin{lemma}{\rm Fix $\mu\in {\bR}$. Suppose $U$ is a finite dimensional  Hilbert space and  $u(t), f(t):[0,L)\ra U$ are two smooth functions   satisfying  the ordinary differential equation
\be
\dot{u}=\mu u +f.
\label{eq: ode}
\ee
Then there exists a constant $C>0$ independent of  $\mu$, $u$ and $f$ such that the following hold.
(a) If $\mu=0$ then
\be
|u(t)-u(t+n)|\leq  \int_t^{t+n}q_{s,L}(f) ds ,\;\;\forall t\in [0, L-n).
\label{eq: a1}
\ee
(b) If   $\mu>0$ then
\be
|u(t)|\leq  e^{-n\mu}|u(t+n)| + \frac{C}{\mu^2}q_{t,L}(f),\;\;\forall t\in [0, L-n).
\label{eq: a2}
\ee
(c) If  $\mu<0$ then}
\be
|u(t+n)| \leq  e^{-n\mu} |u(t)| +\frac{C}{\mu^2}q_{t, L}(f)\,\;\;\forall t\in [0,L-n).
\label{eq: a3}
\ee
\label{lemma: elementary}
\end{lemma}

\noindent {\bf Proof}\hspace{.3cm} We prove only (a) and (b). (c) follows from (b) by time reversal.

\noindent {\sf Proof of (a)} We have
\[
|u(t+1)-u(t)|\leq \int_t^{t_1}|f(s)| ds \leq  \rho_{t}(f).
\]
Thus
\[
|u(t+n)-u(t)|\leq \sum_{k=1}^n|u(t+k)-u(t+k-1)|\leq \sum_{k=1}^n\rho_{t+k-1}(f) \leq q_{t,L}(f).
\]
{\sf Proof of (b)} Denote by ${\bf e}_\mu$ the exponential function $e^{\mu t}$. We have
\[
u(t+1)=e^{\mu}+\int_{0}^1{\bf e}_\mu(1-s)f(t+s) ds
\]
so that by Cauchy's inequality
\[
|u(t+1)-e^{\mu}u(t)|\leq  \rho_0({\bf e}_\mu)\rho_t(f).
\]
Hence
\[
|u(t)|\leq e^{-\mu}|u(t+1)|+e^{-\mu}\rho_0({\bf e}_\mu) \rho_t(f).
\]
Now observe that
\[
e^{-\mu}\rho_0({\bf e}_\mu)= e^{-\mu}\frac{e^{\mu}-1}{\mu} \leq \frac{1}{\mu}.
\]
Set $x_k:= |u(t+k)|$.  The sequence $x_k$ satisfies the  difference inequality
\[
x_k \leq e^{-\mu}x_{k+1} +\frac{1}{\mu}\rho_{t+k}(f).
\]
Thus
\[
x_0 \leq  e^{-n\mu}x_n +\frac{1}{\mu}\sum_{k=0}^{n-1}\rho_{t+k}e^{-(n-1-k)\mu}   \leq  e^{-n\mu} +\frac{1}{\mu}q_{t,L}\sum_{k=0}^{n-1}e^{-k\mu}
\]
\[
\leq e^{-n\mu} +\frac{1}{\mu(1-e^{-\mu})}q_{t,L} \leq e^{-n\mu} +\frac{C}{\mu^2}q_{t,L}(f).
\]
This proves (\ref{eq: a2}) and  the lemma. $\Box$

\bigskip

\noindent {\bf Proof of the Key Estimate}. Let $\hat{u}$ and $\hat{f}$ as in the statement of Proposition \ref{prop: key}.  Using  the spectral decomposition of  $D$  we obtain a  family of  ordinary differential equations of the type (\ref{eq: ode}).  The  Key Estimate is now  an immediate  consequence of Lemma \ref{lemma: elementary}.  The details can be safely left to the reader. $\Box$

\bigskip

\noindent {\bf Proof of Lemma \ref{lemma: asym}}\hspace{.3cm}  Denote by $P_U$ and $P_V$ the orthogonal projections onto $U$ and respectively $V$.

Suppose $\hat{\delta}(U, V) =\sqrt{1-a^2}$, $a \in (0,1)$.  This  means that for every $u \in U$ we have
\[
\|u-P_Vu\|^2\leq (1-a^2)\|u\|^2 
\]
so that
\[
\|P_Vu\|^2 =\|u\|^2 -\|u-P_Vu\|^2 \geq a^2 \|u\|^2
\]
i.e.
\be
\|P_Vu\|\geq  a\|u\|,\;\;\forall u\in U.
\label{eq: low1}
\ee
This shows $P_V$ is one-to-one.

Suppose now that $\hat{\delta}(V, U)= \sqrt{1-b^2}$, $b\in (0,1)$ so that
\[
\delta(U, V) =\max(\sqrt{1-a^2}, \sqrt{1-b^2}) <1.
\]
We deduce similarly
\be
\|P_Uv\|\geq b\|v\|,\;\;\forall v \in V.
\label{eq: low2}
\ee
Introduce the operators
\[
A: U\stackrel{P_V}{\ra}V\stackrel{P_U}{\ra}U,\;\;B:V\stackrel{P_U}{\ra}U\stackrel{P_V}{\ra}V.
\]
Note first  that both $A$ and $B$ are selfadjoint operators.    From (\ref{eq: low1}) and (\ref{eq: low2}) we deduce
\[
\|Au\|\geq ab\|u\|,\;\;\|Bv\|\geq ab\|v\|,\;\;\forall u\in U, \;v\in V.
\]
Thus  both $A$ and $B$ are linear isomorphisms which implies that the operators 
\[
P_V:U\ra V \;\;{\rm and}\;\; P_U: V\ra U
\]
are bounded, one-to-one and onto. We conclude from the closed graph theorem that they must be linear isomorphisms. $\Box$

\begin{remark}{\rm  A similar argument  proves that if $U_r$, $V_r$ is a family of closed subspaces of a Hilbert space $H$  such that
\[
\delta(U_r, V_r)=o(1),\;\;{\rm as}\;\; r\ra \infty
\]
then}
\[
\|{\bf 1}_{U_r} - P_{U_r}P_{V_r}\| + \| {\bf 1}_{V_r} -  P_{V_r}P_{U_r}\|=o(1)\;\;{\rm as}\; \; r\ra \infty.
\]
\label{rem: wae}
\end{remark}

%% file: splice.bbl
\begin{thebibliography}{XXXXX}
 
\addcontentsline{toc}{section}{References}



\bibitem{Andre} A. Andreotti: {\sl Complexes of Partial Differential Operators}, Yale University Press, New Haven and London, 1975.

\bibitem{APS1} M.F. Atiyah,  V.K. Patodi, I.M. Singer: {\sl  Spectral asymmetry and  Riemannian  geometry I},  Math. Proc. Cambridge Philos. Soc. {\bf 77}(1975),  43-69.

\bibitem{BGV} N. Berline, E.Getzler, M. Vergne: {\sl Heat Kernels and Dirac Operators}, Springer Verlag, 1992.

\bibitem{BW} B. Booss-Bavnbek, K. Wojciechovski: {\sl Elliptic Boundary Problems for Dirac Operators}, Birkh\"{a}user,  1993.

\bibitem{BT} R. Bott, L. Tu: {\sl Differential Forms in Algebraic Topology}, Springer Verlag, 1982.


\bibitem{CLM} S. Cappell, R. Lee, E. Miller: {\sl Self-adjoint elliptic operators and manifold decompositions. Part I: Low eigenmodes and stretching}, Comm. Pure Appl. Math., {\bf 49}(1996), 825-866.

\bibitem{Chen} W. Chen:  {\sl A lower bound of the first eigenvalue of certain self-adjoint operators on manifolds containing long necks}, Turkish J. of Math, {\bf 21}(1997), 93-98.


\bibitem{Kato} T. Kato: {\sl  Perturbation Theory for Linear Operators}, Springer Verlag, 1984.


\bibitem{LM} R. Lockart, R. McOwen: {\sl  Elliptic operators on non-compact manifolds}, Ann. Scuola Norm. Sup. Pisa {\bf 12}(1985), 409-446.

\bibitem{Mrow} T. Mrowka: {\sl  A local Mayer-Vietoris principle for Yang-Mills moduli spaces},  PhD Thesis, 1988.

\bibitem{N4} L.I. Nicolaescu: {\sl Rigidity of generalized laplacians and some geometric applications},  Aequationes Mathematic{\ae}, {\bf 48}(1994), 143-162.

\bibitem{N} L.I. Nicolaescu: {\sl  The Maslov index, the spectral flow and decompositions of manifolds}, Duke Math. J., {\bf 80}(1995),  485-533.

\bibitem{N2} L.I. Nicolaescu: {\sl  Generalized Symplectic Geometries and the Index of Families of Elliptic Problems}, Mem. A.M.S., {\bf 128}(1997), no. 609.
 
\bibitem{Yos}   T. Yosida: {\sl  Floer homology and splittings of manifolds}, Ann. of Math. {\bf 134}(1991), 277-324.






\end{thebibliography}
